\documentclass{amsart}

\usepackage[backref=page]{hyperref}       
\usepackage[utf8]{inputenc} 
\usepackage[T1]{fontenc}    
\usepackage{amsfonts}       
\usepackage{nicefrac}       
\usepackage{microtype}      
\usepackage[capitalize]{cleveref}

\hypersetup{
  colorlinks   = true, 
  urlcolor     = teal, 
  linkcolor    = teal, 
  citecolor   = teal, 
  menucolor=teal,%
}

\usepackage[usenames,dvipsnames]{xcolor}

\usepackage[color=teal!70!white, bordercolor=teal, textsize=footnotesize,obeyFinal]{todonotes}

\usepackage{marvosym}
\usepackage{amsmath}
\usepackage{amsgen}
\usepackage{amsbsy}
\usepackage{amsopn}
\usepackage{amsthm}
\usepackage{amscd}
\usepackage{amsfonts}
\usepackage{amssymb}
\usepackage{mathrsfs}

\usepackage{ifsym}

\usepackage{enumerate}
\usepackage{tikz-cd}
\usepackage{epigraph}

\usepackage[usestackEOL]{stackengine}
\renewcommand{\epigraphsize}{\small}
\newcommand{\mytextformat}{\epigraphsize\itshape}
\newcommand{\mysourceformat}{\epigraphsize\scshape}

\let\originalepigraph\epigraph
\renewcommand\epigraph[2]{%
  \setbox0=\hbox{\stackon{\textit{\mytextformat\Longstack{#1}}}%
    {\mysourceformat\scshape\Longstack{#2}}}%
  \ifdim\wd0>.8\linewidth\wd0=.8\linewidth\fi%
  \setlength{\epigraphwidth}{\wd0}%
  \originalepigraph{\textit{#1}}{\textsc{#2}}%
}

\usepackage{color}

\usepackage{graphicx}
\usepackage{comment}

\usepackage{latexsym}

\newtheorem{Theorem}{Theorem}
\newtheorem{Lemma}[Theorem]{Lemma}

\newtheorem{Corollary}[Theorem]{Corollary}

\newtheorem{Proposition}[Theorem]{Proposition}

\newtheorem{Example}[Theorem]{Example}

\def\h{{\rm h}_{\rm top}(g)}
\def\en{{\rm h}_{\rm top}}

\newcommand{\R}{\bf R}
\newcommand{\CC}{{\bf C}}

\newcommand{\cout}[1]{}

\definecolor{lime}{HTML}{A6CE39}
\DeclareRobustCommand{\orcidicon}{%
    \begin{tikzpicture}
    \draw[lime, fill=lime] (0,0)
    circle [radius=0.16]
    node[white] {{\fontfamily{qag}\selectfont \tiny ID}};
    \draw[white, fill=white] (-0.0625,0.095)
    circle [radius=0.007];
    \end{tikzpicture}
    \hspace{-2mm}
}

\begin{document}

\title[Zero entropy on entire Grauert tubes]{Zero entropy on entire Grauert tubes}

\author[P. Su\'arez-Serrato]{P. Su\'arez-Serrato \href{https://orcid.org/0000-0002-1138-0921}{\orcidicon}\\ \today} 

\address{Max Planck Institute for Mathematics, Vivatsgasse 7, Bonn, Germany}

\address{Geometric Intelligence Laboratory, Electrical and Computer Engineering, University of California, Santa Barbara, USA}

\address{{\it and}
Instituto de Matem\'aticas, Universidad Nacional Aut\'onoma de M\'exico UNAM, M\'exico Tenochtitlan}

\begin{abstract}
On a real analytic Riemannian manifold a Grauert tube is an uniquely adapted complex structure defined on the tangent bundle.
It is called entire if it may be defined on the whole tangent bundle.
Here, we show that the geodesic flow of an analytic manifold with entire Grauert tube has zero topological entropy.
Several consequences then follow.
We find that Grauert tubes of convex analytic hypersurfaces of ${\bf R}^3$ are generically finite.
We give a complete classification of $3$-manifolds with entire Grauert tube, showing that these are precisely the $3$-manifolds that admit good complexifications. 
Assuming a manifold has entire Grauert tube, we offer classification statements for several families of smooth $4$-manifolds, determine all such simply connected $5$-manifolds, and offer new topological restrictions for certain manifolds with infinite fundamental groups.
Using these results we present an example of a simply connected and rationally elliptic $5$-manifold which nonetheless does not admit any metric with entire Grauert tube.
\end{abstract}

\maketitle

\section{Introduction and Main Results}

Imagine standing on a manifold with an analytic metric. 
Around each foot, a neighborhood forms, and these overlap as you walk around. 
The transition functions are analytic, so it is reasonable to expect that in your vertical standing direction it is possible to prolong the neighborhood of your standing position to a set so that the augmented coordinate changes there become holomorphic. 
For example, this could be achieved taking power series expansions and making the vertical set small enough to get convergence.  
More formally, a smooth manifold $M$ has a unique real-analytic structure that is compatible with its smooth structure. 
By covering $M$ with charts whose transition functions are real-analytic diffeomorphisms
\[
\phi_{ij}:=\phi_j^{-1}\circ\phi_i , \quad U_{ij}:=\phi_i^{-1}(\phi_i(U_i)\cap\phi_j(U_j))\to U_{ji},
\]
we find open subsets $U_i^{\bf C}\subseteq{\bf C}^n$ with $U_i^{{\bf C}}\cap{\R }^n=U_i$ and $U_{ij}^{{\bf C}}\cap{\R }^n=U_{ij}$ such that the real-analytic $\phi_{ij}$ extend to biholomorphisms $\phi_{ij}^{{\bf C}}:U_{ij}^{{\bf C}}\to U_{ji}^{{\bf C}}$ satisfying the usual cocycle conditions.
Using the cocycle conditions the complexification $M^{{\bf C}}$ is then well defined as a quotient of the disjoint union, $\left(\coprod_i U_i^{{\bf C}}\right)/\sim$.
Here $z_i\sim z_j$ if and only if $z_i\in U_{ij}^{{\bf C}}$ and $z_j = \phi_{ij}^{{\bf C}}(z_i)$.
Coordinate charts $U_i^{{\bf C}}\to M^{{\bf C}}$ with biholomorphic transition functions are induced by the maps $U_i^{{\bf C}}\hookrightarrow\coprod U_i^{{\bf C}}$.
This description is part of the proof of Bruhat--Whitney's theorem showing the existence of $M^{{\bf C}}$ \cite{BruhatWhitney59}.
A similar idea is already present in the 1958 work of Morrey on embeddings of analytic manifolds \cite{Morrey58}.
Such a complexification was also used by Grauert that same year when proving embeddability of abstract real analytic manifolds \cite{Grauert58}.
Even though this local extension of patches paints an intuitive picture of how a complexification may be built, the choices involved in such a construction prevent it from being canonical.

In the early 1990's a canonical construction of such an adapted complex structure, extended within the tangent bundle, was introduced and studied by Sz\H{o}ke \cite{Szoke91}, Guillemin--Stenzel \cite{GuilleminStenzel91}, and Lempert--Sz\H{o}ke \cite{LempertSzoke91}.
They showed that set $T^{r}M$ of tangent vectors of length less than $r$ of a smooth manifold $M$ with a Riemannian metric $g$ admits a K\"ahler structure, for some $r>0$, if and only if the manifold $(M,g)$ is analytic \cite{GuilleminStenzel91, LempertSzoke91}.
These structures are called {\it adapted complex structures} and also {\it Grauert tubes} \cite{Grauert58}.
When such a Grauert tube exists, the adapted complex structure is unique \cite{Lempert93, Szoke91}.
A Grauert tube is called entire if it may be extended to all of $TM$.
For example, over the standard round metric on $S^2$ the complex structure of its Grauert tube
on $TS^{2}$ is biholomorphic to the complex quadric $\{ z^2_1 + z^2_2 + z^2_3 \,=\, 1 \}$ in $\CC^{3}$, and it is entire.

A {\it good complexification} of a closed smooth manifold $M$ is, by definition, a smooth affine algebraic variety $U$ over the real numbers such that $M$ is diffeomorphic to $U(\R)$ and the inclusion $U(\R) \to U({\bf C})$ is a homotopy equivalence \cite{Totaro03}.
Burns conjectured that for every closed Riemannian manifold $M$ with entire Grauert tube, the complex manifold $TM$ is an affine algebraic variety in a natural way \cite{Burns82}. 
Assuming this is correct, then the complex manifold $TM$ would be a good complexification of $M$.
Recall that a smooth affine variety can acquire a symplectic structure through an algebraic embedding into a complex space, followed by the pulling back of the standard symplectic form onto the variety.

In 2003 Totaro showed that all closed manifolds with non-negative sectional curvature known then have good complexifications  \cite{Totaro03}.
Then he asked if, on a given manifold, the existence of a good complexification is equivalent to the existence of a metric of a non-negative sectional curvature. 
This intuition is based on the fact that if a manifold admits an entire Grauert tube, then it admits a smooth metric of nonnegative sectional curvature.
While it would be interesting to check if newer non-negatively curved examples---such as those constructed by Goette--Kerin--Shankar \cite{GoetteKerinShankar20}---also conform to Totaro's equivalence, that quest lies beyond the scope of this paper. 
Moreover, a precise bound on the size of a Grauert tube was provided by Lempert and Sz\H{o}ke, showing that any negative sectional curvature prevents a Grauert tube from covering the entire tangent bundle \cite{LempertSzoke91}.

It appears that all known examples of manifolds with an entire Grauert tube were either already known to be nonnegatively curved, or are constructed by taking quotients of such nonnegatively curved manifolds, as shown by Aguilar \cite{Aguilar07, Aguilar09, Aguilar13}.
Without attempting to be exhaustive, some known results about manifolds with entire Grauert tube include the following.
Burns and Leung showed that the round metric is the only Zoll metric on $S^{n}$ with an entire Grauert tube \cite{BurnsLeung62}.
Patrizio–Wong showed that all compact rank-one symmetric spaces (known as CROSSes) admit entire Grauert tubes \cite{PatrizioWong91}. 
Sz\H{o}ke has even shown that all compact normal Riemannian homogeneous manifolds have entire tubes, and they are quotient spaces of complexified isometry groups \cite{Szoke98}.
Aslam, Burns, and Irvine \cite{AslamBurnsIrvine18}, investigated the existence of entire Grauert tubes on left-invariant metrics $g$ on $SU(2)$, using classical mechanics and integrability properties of the geodesic flow.
They leveraged the complete integrability of the geodesic flow to find an obstruction to tubes being entire.
They showed that if the Grauert tube is entire on $SU(2)$, then the metric $g$ must be homogeneous under the action of $U(2)$.
Moreover, Aslam--Burns--Irvine also showed that if $M$ is a compact connected Lie group $G$ equipped with a left-invariant metric $g$ with an entire Grauert tube, then the tube is biholomorphic to the group complexification $G_{\CC}$ of $G$ \cite{AslamBurnsIrvine18}.
Aguilar \cite{Aguilar13} has found a technical sequence of conditions generalizing previous work of Lempert--Sz\H{o}ke \cite{LempertSzoke91}, producing necessary and sufficient conditions for a Grauert tube to be entire.
These turn out to be, however, hard to compute or interpret in practice.
In a related direction, a recent result of Pali shows integrability properties for geodesic flows of metrics with entire Grauert tube \cite{Pali24}. 

Briefly, a manifold $M$ is called {\it rationally elliptic} if the total rational homotopy $\pi_{*}(M)\otimes {\bf Q}$ is finite dimensional \cite{Paternain99}.
Recently, Chen \cite{Chen24} showed that a simply connected manifold with entire Grauert tube is rationally elliptic.

In this paper, we provide a new dynamical obstruction to the existence of entire Grauert tubes, using the topological entropy of the geodesic flow.
The topological entropy is a rough measure of orbit complexity in a dynamical system.
It quantifies the exponential growth rate of distinct orbits as a function of their length.
Smooth metrics with null entropy often display highly symmetric properties. 
As basic examples, geodesic flows of flat metrics and round metrics have null entropy \cite{Manning79, Paternain99}.

We now state our main result.

\begin{Theorem}\label{thm:entire-Grauert-null-htop}
Let $(M,g)$ be a closed and oriented real analytic manifold with entire Grauert tube.
Then, the topological entropy $h_{top}(g)$ of the geodesic flow of $g$ vanishes.
\end{Theorem}

The proof will follow after we collect various facts and profound properties of topological entropy of geodesic flows.
As a brief summary, the following two ideas are key to proving \cref{thm:entire-Grauert-null-htop}.
In \cref{eqn:Chens-poly-bd} we recover a technical aspect of Chen's proof (initially also used by Aguilar), implying the geodesic counting function of a manifold with entire Grauert tube is bounded above by a polynomial.
We recall Ma\~n\'e's  formula in \cref{eqn:mane} relating the topological entropy of the geodesic flow of a Riemannian metric with the geodesic counting function.
Combining these elements, we prove that a manifold with entire Grauert tube must have zero entropy.

As first observed by Gromov \cite{Gromov87}, \cref{thm:entire-Grauert-null-htop} implies that a simply connected manifold with entire Grauert tube is rationally elliptic (see also \cite[Corollary 5.21]{Paternain99}), recovering Chen's result \cite{Chen24}.

We remark---in relation to \cref{thm:entire-Grauert-null-htop} and the results mentioned above on the integrability of the geodesic flow---that Bolsinov--Taimanov \cite{BolsinovTaimanov00} and Butler \cite{Butler04} have produced and analyzed interesting examples of manifolds with positive entropy integrable geodesic flows. 

By Manning and Dinaburg having null entropy implies that the growth type of the fundamental group is at most subexponential \cite{Manning79, Dinaburg71, Paternain99}.
So we immediately obtain the following Corollary as a consequence:

\begin{Corollary}\label{cor:entire-Grauert-subexp-pi1}
Let $(M,g)$ be a closed and oriented real analytic manifold with entire Grauert tube.
Then the growth of $\pi_1(M)$ is of subexponential type.
\end{Corollary}

The well known fact that the only closed orientable surfaces with entire Grauert tube are $S^2$ and $T^2$ easily follows from \cref{cor:entire-Grauert-subexp-pi1}.
Establishing exactly which metrics on $S^2$ have entire Grauert tube is still an open problem.
Likewise, knowing precisely when a metric on $S^2$ admits a zero entropy metric has been considered a formidable problem \cite{Paternain99}.
Having both perspectives in mind helps us rule out certain positive curvature metrics on ellipsoids from having entire Grauert tube.

\begin{Example}\label{thm:Paternain-metrics-pos-curv-finite-GT}
 
 Let $0<a_1<a_2< a_3$ and consider the ellipsoid $E$ in $\R ^3$ defined by,
 \[
 \frac{x^2}{a_1}+\frac{y^2}{a_2}+\frac{z^2}{a_3}=1.
 \]
 Write $g_{E}$ for the canonical metric of $E\subset \R^3$, and let $r_1, r_2, r_3$ be real numbers.
 Then, for $\varepsilon \neq 0$ sufficiently small, the following analytic Paternain metrics $g_{P}$ on $E$ have positive sectional curvature and positive entropy \cite{Paternain00}:
 \[ 
 g_{P}=\frac{1-\varepsilon(r_1 x + r_2 y + r_3 z)}{a_1 a_2 a_3 \left( \frac{x^2}{a_1^2}+\frac{y^2}{a_2^2}+\frac{z^2}{a_3^2} \right) }
 \]
 By \cref{thm:entire-Grauert-null-htop} the Grauert tubes of the Paternain metrics $g_{P}$ on $E$ are finite.
\end{Example}

\cref{thm:entire-Grauert-null-htop} provides information about the rarity of analytic metrics on $S^2$ with entire Grauert tubes.
A result by Contreras \cite{Contreras10} states that $C^{\infty}$ metrics $g$ generically have geodesic flows with positive topological entropy, $\h>0$.
The case of analytic metrics is somewhat more subtle, as it is not possible to use bump functions and constructions are inherently more delicate.
Nevertheless, Clarke \cite{Clarke22} showed that $C^{\omega}$ metrics $g$ on convex surfaces in ${\bf R}^{3}$ generically have positive $\h$ as well. 
Hence, as a direct consequence of \cref{thm:entire-Grauert-null-htop} we obtain the next result:

\begin{Corollary}\label{cor:Generic-finite-GTs}
Let $(M,g)$ be a convex analytic surface of ${\bf R}^{3}$. Then, $C^{\omega}$-generically, Grauert tubes on $(M,g)$ are finite.
\end{Corollary}

Wielding \cref{thm:entire-Grauert-null-htop}, and standing on the shoulders of Anderson--Paternain \cite{AndersonPaternain03}, Biswas--Mj \cite{BiswasMj15}, Paternain \cite{Paternain93, Paternain00} and Paternain--Petean \cite{PaternainPetean04, PaternainPetean06}, we will harvest the next eight theorems and corollaries.

To begin, we present a classification of $3$-manifolds with entire Grauert tubes in the next result. 
Recall that a model geometry ${\mathbb X}$ is a complete simply connected Riemannian manifold $X$ such that the group of isometries acts transitively on $X$, and it contains a discrete subgroup with compact quotient.
A closed manifold $M$ is said to admit a geometric structure modelled on ${\mathbb X}$ if there exists a Riemannian metric on $M$ such that the Riemannian universal covering of $M$ is $X$.
We write ${\mathbb S}^n$ for the spherical $n$-dimensional geometry of constant positive sectional curvature equal to 1, and ${\mathbb E}^n$ for the flat $n$-dimensional Euclidean geometry.  

\begin{Theorem}\label{thm:3D-EntireGrauert}
 An analytic Riemannian $3$-manifold $(Y,g)$ has entire Grauert tube if, and only if, $(Y,g)$ is a geometric manifold modelled on ${\mathbb S}^3, \, {\mathbb S}^2\times {\mathbb E}$, or ${\mathbb E}^3$.
 Therefore, $(Y,g)$ has entire Grauert tube if, and only if, it has a good complexification.
\end{Theorem}

With \cref{thm:entire-Grauert-null-htop} we can extract information about the $4$-manifolds that admit entire Grauert tubes.
Recall that the Kodaira dimension of a complex manifold $X$ is defined by $$\kappa = \limsup_{m\to \infty} (\log(P_{m}(X)) / \log m),$$ with $P_{m}(X)$ equal to the dimension of the space of holomorphic sections of the $m$-th tensor power of the canonical line bundle of $X$, setting $\kappa := -\infty$ if $P_{m}(X)=0$ for all $m$.
In the following notation we will write $\overline{X}$ to mean the oriented manifold $X$ with the opposite orientation.

\begin{Theorem}\label{thm:4D-simply-connected-EntireGrauert}
Let $(X,g)$ be an analytic Riemannian $4$-manifold with entire Grauert tube.
\begin{enumerate}[(i)]
 \item If $(X,g)$ is simply connected, then $X$ is homeomorphic to $S^4, {\CC}P^2, S^2\times S^2,$ $ {\CC}P^2 \# {\CC}P^2$, or ${\CC}P^2 \# \overline{{\CC}P}^2$. 
\item If $X$ is any known example of a compact complex surface that is not of K\"ahler type, then $X$ is diffeomorphic to a Hopf surface, and it admits a geometric structure modeled on ${\mathbb S}^3\times {\mathbb E}^3$.
\item If $X$ is a compact complex K\"ahler surface with Kodaira dimension at most 1, then $X$ admits a geometric structure modelled on ${\CC}P^2, {\mathbb S}^2\times {\mathbb S}^2, {\mathbb S}^2\times {\mathbb E}^2$, or ${\mathbb E}^4$.
Moreover, the Kodaira dimension of $X$ is equal to either $-\infty$ or $0$, and $X$ is diffeomorphic to one of: ${\CC}P^2$, a ruled surface of genus $0$ or $1$, a complex torus, or a hyperelliptic surface.
\end{enumerate}
\end{Theorem}

McLean studied when a (co)tangent space is symplectomorphic to an affine variety using the growth rate of wrapped Floer cohomology \cite{McLean18}.
He showed that if the (co)tangent bundle of a simply connected 4-manifold is symplectomorphic to a smooth affine variety, then it is homeomorphic to one in the list given in item $(i)$ of \cref{thm:4D-simply-connected-EntireGrauert}.
This result is significant in relation to Burn's conjecture, and shows that Burn's conjecture remains compatible with our results here as well.

 A complex surface $S$ is an elliptic surface if there exists a holomorphic map $S\to C$ to a complex curve $C$ with general fiber an elliptic curve. 
 A good orbifold is one constructed as a global quotient of a smooth manifold by the action of a discrete group.
We next characterize when a compact complex surface $S$ with an elliptic fibration $\pi: S\to C$ over a good orbifold $C$ has entire Grauert tube.

\begin{Theorem}\label{thm:4D-elliptic-surfaces-EntireGrauert}
Let $S$ be an elliptic surface with a fibration $\pi: S\to C$, over a good orbifold $C$. 
Then $S$ has entire Grauert tube if, and only if, it is one of the following:
\begin{enumerate}[(i)]
    \item A geometric manifold modelled on ${\mathbb S}^3\times {\mathbb E}$, in which case $S$ is a Hopf surface.
    \item A geometric manifold modelled on ${\mathbb S}^2\times {\mathbb E}^2$, in which case $S$ is a ruled surface of genus $1$.
    \item A geometric manifold modelled on ${\mathbb E}^4$, in which case $S$ is a torus or a hyperelliptic surface.
\end{enumerate}
\end{Theorem}

In the case of $4$-manifolds with infinite fundamental group the following result holds true.

\begin{Theorem}\label{thm:infinite-pi1-EntireGrauert}
Let $X$ be a closed analytic $4$-manifold with entire Grauert tube and infinite fundamental group of polynomial growth type.
Then $X$ is finitely covered by a manifold homeomorphic to $S^3\times S^1$, $S^2\times T^2$, or it is a flat manifold.
\end{Theorem}

Smale \cite{Smale62} and Barden \cite{Barden65}  classified simply connected smooth 5-manifolds, showing they are uniquely characterized by their second homology group and an invariant related to their Stiefel-Whitney class.
Relying on the Barden--Smale classification of smooth simply connected $5$-manifolds, we offer the following result.

\begin{Theorem}\label{thm:simply-conn-5mfds-EGT}
 A closed smooth analytic simply connected orientable $5$-manifold $(M,g)$ has entire Grauert tube if, and only if, it  is diffeomorphic to one of $S^5$, $S^3\times S^2$, the nontrivial $S^3$-bundle over $S^2$, or the Wu manifold ${\rm SU}(3)/{\rm SO}(3)$.
 Therefore, $(M,g)$ has entire Grauert tube if, and only if, it has a good complexification.
\end{Theorem}

Observe that the manifolds featured in \cref{thm:simply-conn-5mfds-EGT} are all the known examples of simply connected $5$-manifolds which admit metrics of non-negative sectional curvature.
Returning to Burn's conjecture, McLean \cite{McLean18} also showed that if the (co)tangent bundle of a simply connected $5$-manifold  $M$ is symplectomorphic to a smooth affine variety, then $M$ is diffeomorphic to one of the manifolds in \cref{thm:simply-conn-5mfds-EGT}.

In the following result we find topological restrictions to a $5$-manifold with infinite fundamental group and entire Grauert tube.

\begin{Corollary}\label{thm:pi1-infinite-5mfds-EGT}
Let $(M,g)$ be a closed analytic $5$-manifold with infinite fundamental group and entire Grauert tube.
Then, the rational homotopy class of the universal covering $\widetilde{M}$ is that of a $1$-connected elliptic complex.
\end{Corollary}

\begin{proof}
Paternain and Petean \cite[Theorem C]{PaternainPetean06} have shown this is the case, once we know $\h =0$, which follows by Theorem \ref{thm:entire-Grauert-null-htop}.
\end{proof}

Next we present a new topological condition that limits the existence of entire Grauert tubes. 
This constraint is particularly useful for manifolds with non-trivial fundamental groups, providing insights into how manifolds with entire Grauert tubes may decompose into a connected sum.

\begin{Corollary}\label{thm:connected-sums-EntireGrauert}
Let $M$ be a closed manifold with  $\dim(M) > 3$.
Suppose that $M$ can be decomposed as a connected sum $M = X_1 \# X_2$, where the order of the $\pi_1(X_1)$ is greater than 2.
If the Grauert tube of $M$ is entire then $X_2$ is a homotopy sphere.
\end{Corollary}

\begin{proof}  
 This follows directly from \cref{thm:entire-Grauert-null-htop}, and a result of Paternain and Petean \cite[Theorem D]{PaternainPetean04}, which reaches the same conclusion for such connected sums, assuming the existence of a metric $g$ with $\h =0$.
\end{proof}

 Let $\Omega M$ denote the loop space of the manifold $M$, that is, the space of continuous maps from $S^1$ to $M$.
The homology $H_{\ast}(\Omega M, {\bf Q})$ of $\Omega M$ is a well defined object \cite{Paternain99}.
Set $b_{j}(\Omega M, {\bf Q}):= \dim H_{j}(\Omega M, {\bf Q})$, and 
 \[
 \mu_{i}(M) := \sum\limits_{j\leq i} b_{j}(\Omega M, {\bf Q}).
 \]
For a closed simply conected manifold $M$, the series $\{ \mu_{i}(M) \}$ 
may grow exponentially or polynomially in $i$.
We say that the growth of the homology of the loop space $H_{\ast}(\Omega M, {\bf Q})$ of $M$ is exponential or polynomial, depending on how $\{ \mu_{i}(M) \}$ grows.
Paternain and Petean found that the type of growth of $H_{\ast}(\Omega M, {\bf Q})$ imposes restrictions to the existence of a null entropy metric \cite{Paternain99, PaternainPetean03}, and hence implies the next result.

\begin{Corollary}\label{thm:LoopSpaceHomology-EntireGrauert}
    Let $(M,g)$ be a closed simply connected and analytic manifold with entire Grauert tube. Then $H_{\ast}(\Omega M, {\bf Q})$ grows polynomially. 
\end{Corollary}

\begin{proof} 
By yet another result of Paternain--Petean \cite[Theorem 8.3]{PaternainPetean03},  if the loop space homology of a manifold $M$ grows exponentially with respect to any field of coefficients, then $M$ does not admit a zero entropy metric.
The result then follows from \cref{thm:entire-Grauert-null-htop}.
\end{proof}

Once more, in relation to Burn's conjecture, MacLean showed that the (co)tangent bundle of a manifold whose free loop space homology grows exponentially is not symplectomorphic to an affine variety \cite{McLean12}. 
So that Burn's conjecture remains compatible with \cref{thm:LoopSpaceHomology-EntireGrauert}.
McLean also reformulated Totaro's question in terms of exponential growth, asking if having  exponential growth prevents a simply connected manifold from admitting a good complexification.
We contribute to this affirming that it if would then, by \cref{thm:LoopSpaceHomology-EntireGrauert}, such a good complexification could not be an entire Grauert tube.

A connected $CW$ complex is said to be nilpotent if its fundamental group is a nilpotent group and it acts nilpotently on all the higher homotopy groups.
Paternain and Petean \cite[Corollary]{PaternainPetean06} showed that if a closed manifold $M$ is a nilpotent $CW$ complex and has a zero entropy metric, then $\pi_{\ast}(\Omega M)\otimes {\bf Q}$ is finite dimensional. 
Consequently, from \cref{thm:entire-Grauert-null-htop} we obtain the next result:

\begin{Corollary}\label{cor:entire-Grauert-nilpotent-CW-null-htop}
Let $M$ be a closed analytic manifold $M$ that is a nilpotent $CW$ complex and has entire Grauert tube.
Then $\pi_{\ast}(\Omega M)\otimes {\bf Q}$ is finite dimensional.
\end{Corollary}

By comparison, Chen's result on the rational ellipticity of manifolds with entire Grauert tubes \cite{Chen24} is restricted to simply connected manifolds. 
This limitation stems from the reliance on arguments about the topology of the loop space, which could be at most extended to manifolds with finite fundamental groups.

Finally, using \cref{thm:LoopSpaceHomology-EntireGrauert}, we present an example of a manifold that is rationally elliptic, and yet its Grauert tubes are never entire. 

\begin{Example}\label{ex:M_2}
    Consider the $5$-dimensional Brieskorn variety $M_2$, defined by the intersection of $S^7\subset {\bf C}^4$ with the zero set of $z_1^{2} +z_2^{3} +z_3^{3} +z_4^{3}$.
    It is a simply connected spin manifold with second homology group isomorphic to ${\bf Z}_2 \oplus {\bf Z}_2$.
    Its loop space homology grows exponentially with ${\bf Z}_2$ coefficients. 
    Therefore, by \cref{thm:LoopSpaceHomology-EntireGrauert} it does not admit an analytic metric with entire Grauert tube. 
    As pointed out by Paternain and Petean \cite{PaternainPetean03}, the manifold $M_2$ has the rational cohomology ring of $S^5$. 
    So the loop space homology of $M_2$ with rational coefficients is bounded, and thus $M_2$ is rationally elliptic.
\end{Example}

A manifold of cohomogeneity-one admits a smooth Lie group action with one-dimensional orbit space.
Cohomogeneity-one actions provide a systematic approach to constructing nonnegatively curved metrics on manifolds by exploiting the symmetries of the group action.
 For example, Hoelscher classified compact simply connected cohomogeneity-one $5$-manifolds \cite[Theorem C]{Hoelscher10}, thus proving that $M_2$ is not of cohomogeneity-one. 
The question of whether cohomogeneity-one manifolds must always admit a null entropy metric remains an intriguing one. 
Do cohomogeneity-one manifolds have entire Grauert tube? 


\cref{sec:context} sets the context and includes preliminary information. \cref{sec:proofs} contains the proofs of the statements above.
\medskip

{\bf Acknowledgements.} I wish to warmly thank the Max Planck Institute for Mathematics in Bonn for providing an excellent working environment and for support during the summer of 2024 when this work was conceived and mostly written up. 
I would also like to thank Werner Ballmann for conversations regarding these topics, and Gabriel Paternain and Jimmy Petean for correspondence about $M_2$, and for their support throughout my career.
I am grateful to Burt Totaro for correspondence explaining the biquotient presentation of the nontrivial $S^3$-bundle over $S^2$, and to Martin Kerin for pointing out relevant references and catching some typos, all their comments helped improve the first version.

\section{Context and preliminaries}\label{sec:context}

Next we will review the definitions and properties of tangent bundles and geodesic counting functions needed in the proofs of our results.
Let $(M,g)$ be a compact smooth Riemannian manifold.

\subsection{Vertical and horizontal subbundles}
Here we recall useful facts about the geometry of the tangent bundle $TM$, closely following Paternain's book \cite{Paternain99}.
Let $\pi: TM \rightarrow M$ be the canonical projection map, and $i: T_x M \rightarrow TM $ be the inclusion.
Set $\pi: TM \rightarrow M$ to be the canonical projection, i.e., $\theta=(x,v)$ in $TM$ implies $\pi(\theta)=x$. 
The canonical vertical subbbundle of $TTM$, with fiber at $\theta$ given by the tangent vectors of curves $\sigma: (-\epsilon, \epsilon) \rightarrow TM$ of the form: $\sigma(t)=(x, v+t \omega)$, with $\omega \in T_x M$, so that $V(\theta)=\ker ((\pi_*)_{\theta})$.
Fix $\xi \in T_{\theta} TM$, and let $z: (-\epsilon, \epsilon)\rightarrow TM$ be an adapted curve to $\xi$, with initial conditions $z(0)=\theta$, and  ${z}'(0)=\xi$.
This data gives rise to a curve $\alpha: (-\epsilon, \epsilon) \rightarrow M, \alpha:=\pi \circ z$, and a vector field $Z$ along $\alpha$, equivalently, $z(t)=(\alpha(t), Z(t)).$

Consider the linear isomorphism $P_t: T_x M \rightarrow T_{\alpha(t)} M$ defined by parallel transport along $\alpha$.
Define the connection map $K: TTM \rightarrow TM$, as,
\[
K_{\theta}(\xi):=(\nabla_{\alpha}Z) (0)=\lim_{t \rightarrow 0} \frac{(P_t)^{-1}Z(t)-Z(0)}{t}.\]

The horizontal subbundle of $TTM$ is the subbundle whose fiber at $\theta$ is given by $H(\theta)=\ker K_{\theta}$, and which can be equivalently constructed using the the horizontal lift
$L_{\theta}: T_x M \rightarrow T_{\theta} TM,$ defined as $\theta=(x, v)$.
Suppose $\omega \in T_x M$, and let $\alpha: (-\epsilon, \epsilon) \rightarrow M$ be such that
$\alpha(0)=x, {\alpha}'(0)=\omega$.
In this case we say that $\alpha$ is adapted to $\omega$.
Let $Z(t)$ be the parallel transport of $v$ along $\alpha$ and $\sigma: (-\epsilon, \epsilon) \rightarrow TM$
be the curve $\sigma(t)=(\alpha(t), Z(t))$. Then define
$$L_{\theta}(w)={\sigma}'(0)\in T_{\theta}TM.$$

\begin{Proposition}[See \cite{Paternain99}]\label{prop:jac}
The maps $K_{\theta}$ and $L_{\theta}$ have the following properties,
$$ (\pi_*)_{\theta} \circ L_{\theta}=Id, \quad {\rm and} \quad K_{\theta} \circ i_*=Id.$$
 Moreover, $T_{\theta}TM= H(\theta)\oplus V(\theta)$,
and the map $j_{\theta}: T_{\theta}TM  \rightarrow T_x M \times T_x M$ given by
$j_{\theta} (\xi)=( (\pi_*)_{\theta} (\xi), K_{\theta}(\xi))$ is a linear isomorphism.
\end{Proposition}

For each $\theta \in TM$, there exists a unique geodesic $\gamma_{\theta}$ in $M$ with initial condition $\theta$.
Let $\xi \in T_{\theta} TM$, and $z: (-\epsilon, \epsilon)\rightarrow TM$ be an adapted curve to $\xi$, (such that $z(0)=\theta,$ and ${z}'(0)=\xi$).
Call $\phi_t$ the geodesic flow of $TM$, then the map
$(s, t) \mapsto \pi \circ \phi_t (z(s))$ gives rise to a variation of $\gamma_{\theta}$.
 The curves $t \mapsto \pi \circ \phi_t (z(s))$ are geodesics,
and therefore the corresponding variational vector fields
$$J_{\xi}:=\frac{\partial}{\partial s}|_{s=0} \pi \circ \phi_t (z(s))$$
are Jacobi fields with initial conditions $J_{\xi}(0)=(\pi_*)_{\theta} (\xi)$, and ${J}'_{\xi}(0)=K_{\theta}({\xi})$.


\subsection{Complex structures on $TM$}
Here we will explain the complex structures that are adapted to $TM$.
We shall begin by describinga foliation on $TM\setminus M$ by Riemannian surfaces.
 Denote by $N_{\tau}: TM \rightarrow TM $ the smooth mapping defined by multiplication by a real $\tau$ on the tangent fibers.
 Let $\gamma: \R \rightarrow M$ be a geodesic. Then, define an immersion $\phi_{\gamma}: \CC \rightarrow TM$ by
 $\phi_{\gamma}(\sigma+ i \tau)=N_{\tau}{\gamma}'(\sigma)$.
 
 Given a couple of geodesics $\gamma$ and $\delta$, if the images $\phi_{\gamma}({\bf C}\setminus {\bf R})$ and $\phi_{\delta}({\bf C}\setminus {\bf R})$
 intersect, then $\gamma$ and is equal to $\delta$, possibly with different velocities. 
 This implies $\phi_{\gamma}({\bf C})=\phi_{\delta}({\bf C}).$
 Hence the images of ${\bf C}\setminus {\bf R}$ under  $\phi_{\gamma}$ define a smooth surface foliation of $TM \setminus M$. 
 Notice that each leaf carries a complex structure, inherited from ${\bf C}$ through $\phi_{\gamma}$. 
 As we extend the leaves to a foliation of $M$, this extended foliation may become singular. 
 
Let $r>0$ and define $T^r M=\{{v\in TM \, | \, g(v,v)<r^2}\}$.
We will say that a smooth complex structure on $T^r M$ is {\bf adapted} if the leaves of the foliation $\mathcal{F}$ with the complex structure inherited from $\CC$ are complex submanifolds of $T^r M$.
 The following Theorem was shown by Guillemin and Stenzel \cite{GuilleminStenzel91}, Lempert and Sz\H{o}ke \cite{LempertSzoke91}, and Sz\H{o}ke\cite{Szoke91}.

\begin{Theorem} Let $M$ be a compact real analytic manifold equipped with a real analytic metric $g$.
  Then there exists some $r>0$ such that $T^r M$ carries an unique adapted complex structure.
 \end{Theorem}

 The manifold $M$ is said to have entire  Grauert  tube when the adapted complex
 structure is defined on the whole tangent bundle, that is $r=\infty$.
 
 We will now describe the adapted complex structure on $T^r M$.
 Let $\theta \in T^r M \setminus M$, and $x=\pi(\theta)$, and $\gamma$ be a geodesic determined by $\theta$.
 Choose tangent vectors $v_1, v_2, \cdots, v_{n-1}$ such that the set $v_1, v_2,  \cdots, v_{n-1}, v_n:=\frac{{\gamma}'(0)}{|{\gamma}'(0)|}$ is an orthonormal basis of $T_x M$.
Write $L_{\theta}$ to denote the leaf of the foliation $\mathcal{F}$ over $\theta$.
 A vector $\bar{\xi} \in T_{\theta} TM$ determines the parallel vector field $\xi$ along $L_{\theta}$, and by definition it is invariant under both the actions of the semi-group of $N_{\tau}$ and of the geodesic flow $\varphi_{t}$ of $g$.
 Therefore, ${\xi}|_{\R}$ is a Jacobi field along $\gamma$.
Choose a set of vectors $\bar{\xi}_1, \bar{\xi}_2, \cdots, \bar{\xi}_{n}, \bar{\eta}_1, \bar{\eta}_2, \cdots, \bar{\eta}_{n} \in T_{\theta} TM$ satisfying
 $$ (\pi_*)_{\theta} (\bar {\xi}_j)=v_j,\, \ K_{\theta} (\bar {\xi}_j)=0,\, (\pi_*)_{\theta} (\bar {\eta}_j)=0, \, {\rm and}\, \ K_{\theta} (\bar {\eta}_j)=v_j.$$
Using  the connection map $K: TTM \rightarrow TM$ described in section $2$, extend $\bar {\xi}_j$ and $\bar {\eta}_j$, obtaining parallel vector fields
 $\xi_1, \xi_2, \cdots, \xi_n, \eta_1, \eta_2, \cdots, \eta_n$ along $L_{\theta}$.
In this way the Jacobi fields $\xi_1|_{\bf R}, \xi_2|_{\bf R}, \cdots, \xi_n|_{\bf R}$ are  linearly independent, away from a discrete subset $S_1$ of ${\bf R}$.
 Therefore, there exist smooth real valued functions $\phi_{jk}$ defined on ${\bf R} \setminus S_1$ such that
 \[
 \eta_k|_{{\bf R}}=\sum_{j=1}^n \phi_{jk} \xi_j|_{{\bf R}}.
 \]
 Relying on the adapted complex structure, the functions $\phi_{jk}$ admit meromorphic extensions $f_{jk}$ over the domain
  \[
 D= \left\{ {\sigma + i \tau \in {\bf C}\, : \,  \ |\tau|< \frac{r}{\sqrt{g(\theta, \theta)}}} \right\}
 \]
 such that for each $j, k$, the poles of $f_{jk}$ lies on ${\bf R}$ and the matrix ${\rm Im} (f_{jk})|_{D\setminus {\bf R}}$ is invertible.
 Let $(e_{jk})=({\rm Im} f_{jk}(i))^{-1}$. Therefore, the complex structure $J$ satisfies
 $$J \bar{\xi}_h=\sum_{k=1}^{n} e_{kh}\times
  \left[ \bar{\eta}_k-\sum_{j=1}^{n} {\rm Re}  f_{jk}(i) \bar{\xi}_j\right].$$

  Observe that the fields 
  \[
  \xi_1|_{\bf R}, \xi_2|_{\bf R}, \cdots, \xi_{n-1}|_{\bf R}, \eta_1|_{\bf R}, \eta_2|_{{\bf R}}, \cdots, \eta_{n-1}|_{\bf R}
  \]
  are normal Jacobi fields, and that $\xi_n|_{\bf R}, \eta_n|_{\bf R}$ are tangential Jacobi fields. 
  Therefore, for $ \ 1 \leq j, k \leq n-1$,
  \[
  \phi_{nk}=\phi_{jn}\equiv 0, \quad f_{nk}=f_{jn}\equiv 0,  \quad {\rm and} \quad e_{nk}=e_{jn}\equiv 0.\]

 Consider the $n$-tuples
 $$\Xi=(\xi_1, \xi_2, \cdots, \xi_n), \ H=(\eta_1, \eta_2, \cdots, \eta_n)$$
 and holomorphic $n$-tuples
 $$\Xi^{1,0}=(\xi_1^{1,0}, \xi_2^{1,0}, \cdots, \xi_n^{1,0}),  \ H^{1,0}=(\eta_1^{1,0}, \eta_2^{1,0}, \cdots, \eta_n^{1,0}),$$
 where $\xi_j^{1,0}=\frac{1}{2}(\xi_j- i J \xi_j)$ and $J$ is the adapted complex structure.

Then we have (see \cite{Szoke91, Chen24}),
 $$H(\sigma)=\Xi(\sigma) f(\sigma),$$
 $$H^{1,0}(\sigma+ i \tau)=\Xi^{1,0} (\sigma+ i \tau) f(\sigma + i \tau)\,\, {\rm and},$$
  $$f(\sigma+ i \tau)=(f_{jk}(\sigma+ i \tau)), \ \sigma \in {\bf R}\setminus S_1, \ |\tau| < \frac{r}{\sqrt{g(\theta, \theta)}}.$$

The following results are proved by Lempert--Sz\H{o}ke  \cite{LempertSzoke91}, and by Sz\H{o}ke \cite{Szoke91}.

\begin{Proposition} \label{prop:linear}
(1) The vectors $\xi_1^{1,0}, \xi_2^{1,0}, \cdots \xi_n^{1,0}$ are linearly independent over ${\bf C}$ on $D \setminus {\bf R}$.
The same is true for the vectors $\eta_1^{1,0}, \eta_2^{1,0}, \cdots, \eta_n^{1,0}$.
\\
(2) The $2n$ vectors $\xi_j, \eta_k$ are linearly independent in points $\sigma+ i \tau \in D \setminus {\bf R}$.
\end{Proposition}

\begin{Theorem} \label{thm:mero}
The matrix valued meromorphic function $f(\sigma + i \tau )$ is symmetric (as a matrix) and satisfies  $f(0)=0, f'(0)={\rm Id} $.
Moreover, if $\sigma + i \tau \in D, \tau >0$, then ${\rm Im} f(\sigma + i \tau)$ is a symmetric, positive definite matrix.
\end{Theorem}

\subsection{Growth rate of counting functions}

Let $x$ be a point in $M$, and $T >0$.
Define 
\[
D_T :=\{{ v \in T_x M \, | \, g(v,v)\leq T^2}\}\]
 to be the disk of radius $T$ in $T_x M$. 
 Define the counting function $n_T (x, y)$ by
 \[
 n_T (x, y):=\sharp \left( (\exp_x)^{-1}(y)\cap D_T \right).\]
 
 Observe that $n_T (x,y)$ counts the number of geodesic arcs of length $\leq T$ joining $x$ to $y$.

\par
 The following result, summarizing the contributions of Berger--Bott \cite{BergerBott62}, Gromov \cite{Gromov78}, and Paternain \cite{Paternain99}, is crucial for the proof of our main Theorem.
\begin{Theorem} [\cite{BergerBott62, Gromov78,Paternain99}] \label{thm:Gromov-count-fn}
\begin{equation} \label{eqn:gromov-paternain}
\int_M n_T (x, y)\, dy =\int_0^T  d \sigma \int_{S_{x}} \sqrt{\det (g(J_j (\sigma), J_k (\sigma)))_{j, k=1,2, \cdots, n-1}} \ d\theta.
\end{equation}
Here $S_{x}$ is the unit sphere of $T_x M$.  Moreover,
$J_j, j=1, 2, \cdots, n-1$ are Jacobi fields along the unique geodesic $\gamma$ determined by $\theta$ in $U$ (i.e. $\gamma(0)=x, {\gamma}'(0)=\theta$) with initial conditions, $J_j(0)=0$ and  ${J}'_j (0)=v_j$,
where $v_j, j=1, 2 \cdots n-1$ is an orthonormal basis of $T_{\theta} U$.
\end{Theorem}

Chen showed that when $M$ has  entire Grauert tube the right hand side in \cref{eqn:gromov-paternain} can be further described by
a  matrix valued holomorphic function on the upper half plane \cite{Chen24}. 
Applying Fatou's representation theorem to this function, and relying on \cref{prop:linear} and \cref{thm:mero} Chen found that $\int_M n_T (x, y) dy$ has  polynomial growth and then used this to conclude that $M$ is topologically elliptic.

We include the following result in the form we need, the arguments that prove it may be consulted in Chen's work. 

\begin{Proposition}\label{prop:Chens-polynomial-bound}
Let $M$ be a closed $n$-manifold with entire Grauert tube. Then there exists a polynomial $P(T)$, of degree at most $n$, such that 
\begin{equation}\label{eqn:Chens-poly-bd}
\int_M n_T (x, y)\, dy \leq P(T).
\end{equation}
\end{Proposition}

\begin{proof}
    See the proof of \cite[Theorem 1.1]{Chen24}.
\end{proof}

\subsection{Geodesic Flows.}

We refer interested readers to G.P. Paternain's book \cite{Paternain99} on geodesic flows for a comprehensive treatment of the concepts in this section.

Let $R$ denote the Riemann tensor of the Riemannian manifold $(M,g)$.
 A vector field $J$ along a geodesic $\gamma: (a,b)\rightarrow M$ is called a {\emph{Jacobi field}} if it satisfies
$$\ddot{J} + R(\dot{\gamma}, J)\dot{\gamma}=0.$$

We can describe the differential of the geodesic flow in terms of Jacobi fields. 
As before, we denote by $\pi:TM\rightarrow M$ the canonical projection and by $K:TTM \rightarrow M$ the connection map.
 Let $\xi \in T_{\theta}TM$ and $z:(-\varepsilon , \varepsilon ) \rightarrow TM$ be an adapted curve to $\xi$.
The curves $t\mapsto \pi\circ \phi_{t}(z(s))$ are geodesics and therefore the vector field
$$J_{\xi}(t) :=  \frac{\partial}{\partial s}\mid_{s=0}(\pi\circ \phi_{t}(z(s)))$$
 is a Jacobi vector field.

\begin{Lemma}\label{lem:difgd} {{\rm (\cite[p.20]{Paternain99})} Given $\theta \in  TM$, $\xi \in T_{\theta}TM$ and $t\in \R$, we have
$$d_{\theta}\phi_{t}(\xi ) = (J_{\xi}(t), \dot{J}_{\xi}(t)).$$ }
\end{Lemma}

Observe that, informally, \cref{lem:difgd} means that Jacobi fields determine the growth of geodesics.


\subsection{Topological Entropy.}
Here we recall the definition of topological entropy of a dynamical system.
We recommend the comprehensive book by Hasselblatt and Katok \cite{HasselblattKatok95} for further exploration of these dynamical ideas.
The topological entropy quantifies the exponential growth rate of the number of essentially different orbit segments of a certain length.
It measures the complexity of the orbit structure of a dynamical system.

Let $(X,d)$ be a compact metric space and $f:X\rightarrow X$ a continuous map. 
For each $n$ in $\mathbf{N}$, 
\[
d_{n}(x,y)=\max\limits_{0\leq k \leq n-1} d(f^{k}(x), f^{k}(y)  )\]
measures the maximum distance between iterates of $x$ and $y$.
Each $d_{n}$ is a metric on $X$, $d_{n}\geq d_{n-1}$ and they all induce the same topology on $X$.
Fix $\epsilon >0$, a subset $A\subset X$ is $(n,\epsilon)$-spanning if for every $x\in X$ there is $y\in A$ such that $d_{n}(x,y)<\epsilon$. 
By compactness, there are a finite number of elements in an $(n,\epsilon)$-spanning set. 
Let ${\rm{span}}(n, \epsilon , f)$ be the minimum cardinality of an $(n,\epsilon)$-spanning set.
A subset $A\subset X$ is $(n,\epsilon)$-separated if any two distinct points in $A$ are at least $\epsilon$ apart in the metric $d_{n}$. 
Any $(n,\epsilon)$-separating set is finite. Let ${\rm{sep}}(n, \epsilon , f)$  be the  maximum cardinality of an $(n,\epsilon)$-separated set.
Let ${\rm{cov}}(n, \epsilon , f)$ be the minimum cardinality of a covering of $X$ by sets of $d_{n}$-diameter less than $\epsilon$. Again by compactness, ${\rm{cov}}(n, \epsilon , f)$  is finite.
The values of ${\rm{span}}(n, \epsilon , f)$, ${\rm{sep}}(n, \epsilon , f)$ and ${\rm{cov}}(n, \epsilon , f)$  count the number of orbit segments of length $n$ that are distinguishable at scale $\epsilon$. They are related in the following way, 
\begin{equation}\label{eqn:span-sep-cov-ineq}
  {\rm{cov}}(n, 2\epsilon , f) \leq {\rm{span}}(n, \epsilon , f) \leq {\rm{sep}}(n, \epsilon , f) \leq {\rm{cov}}(n, \epsilon , f).
\end{equation}

Let
$$h_{\epsilon}(f)=\limsup\limits_{n\rightarrow \infty} \frac{1}{n}\log ( {\rm{cov}}(n, \epsilon , f) ).$$

The quantity ${\rm{cov}}(n, \epsilon , f)$ increases monotonically as $\epsilon$ decreases, so $h_{\epsilon}(f)$ does so as well.
Thus the limit
$$\en(f)=\lim\limits_{\epsilon \rightarrow 0^{+} }h_{\epsilon}(f)$$
exists, and it is called the {\emph{topological entropy} of $f$.

The inequalities in \cref{eqn:span-sep-cov-ineq} above imply that  equivalent definitions can be given using ${\rm{span}}(n, \epsilon , f)$ or ${\rm{sep}}(n, \epsilon , f)$.
\begin{eqnarray*}
\en(f) & = &\lim\limits_{\epsilon \rightarrow 0^{+} } \limsup\limits_{n\rightarrow \infty} \frac{1}{n}\log ( {\rm{span}}(n, \epsilon , f)) \\
& = & \lim\limits_{\epsilon \rightarrow 0^{+} } \limsup\limits_{ n \rightarrow \infty}  \frac{1}{n}\log ({\rm{sep}}(n, \epsilon , f) ).
\end{eqnarray*}

The topological entropy of a continuous flow can be defined as the entropy of the time one map. Alternatively, it can also be defined using the analog $d_{T}$, $T>0$ of the metrics $d_{n}$. The definitions are equivalent because of the equicontinuity of the family of time $t$-maps, $t\in [0,1]$.

\subsection{Topological Entropy of the Geodesic Flow.}

Recall that the geodesic flow $\varphi_{t}$ leaves the unit sphere bundle $SM$ invariant. 
When $M$ is a compact manifold, then $SM$ is also compact.
We compute the topological entropy $\h(\varphi_{t}|_{SM})$ of  $\varphi_{t}$ restricted to $SM$, and call this value the topological entropy of the geodesic flow \cite{Manning79}. 
%
Let $d$ be any distance function
compatible with the topology of $SM$. For each $T>0$ we define a new
distance function,
$$d_{T}(x,y):= \max\limits_{0\leq t\leq T} d(\phi_{t}(x),\phi_{t}(y)).$$

As $SM$ is compact, we can consider the minimal number $N(\epsilon, T)$ of balls of radius $\epsilon >0$ in the metric $d_{T}$ that are necessary to cover $SM$, and define,
$$h(\phi, \epsilon):= \limsup\limits_{T\rightarrow
\infty}\frac{1}{T}\log N(\epsilon, \phi).$$

Notice that the function $\epsilon\mapsto h(\phi, \epsilon)$ is monotone decreasing as $\epsilon\mapsto 0$, and therefore the following limit exists,
$$\h:= \lim\limits_{\epsilon\rightarrow 0}h(\phi, \epsilon).$$

The number $\h$ is the {\emph{topological entropy}} of the geodesic flow $\varphi_{t}$ of $g$.
The following formula, due to R. Ma\~{n}\'{e}, gives $\h$ a Riemannian metric a clear geometric meaning. 
\begin{Theorem} {\rm {(Ma\~{n}\'{e}'s Formula \cite{Paternain99}) }} Given points $x$ and $y$ in $M$ and $T>0$, define $n_{T}(x,y)$ to be
the number of geodesic arcs joining $x$ and $y$ with length $\leq T$,
then:
\begin{equation}\label{eqn:mane}
\h=\lim\limits_{T \rightarrow \infty} \frac{1}{T}\log \int\limits_{M\times M}n_{T}(x,y)
\, dx\,  dy
\end{equation}
\end{Theorem}

Another useful equivalence is found in the following formula \cite[p.83]{Paternain99}: 

\begin{equation}\label{eqn:htop-Jacobi-formula}
\h = \lim\limits_{T \rightarrow \infty} \frac{1}{T}\log \int\limits_{SM} | \det ( d_{\theta}\phi_{T}|_{V(\theta )} ) | \, d\theta
\end{equation}
Here, $V(\theta )$ is the tangent space to the fiber $T_{x}M \subset TM$ at the point $\theta $.
From \cref{eqn:htop-Jacobi-formula} and lemma \ref{lem:difgd}, we see that $\h =0$ when the Jacobi fields of a metric $g$ grow at most polynomially.

%
%

\section{Proofs}\label{sec:proofs}

\begin{proof}[Proof of \cref{thm:entire-Grauert-null-htop}.]

Having explained the growth of geodesics of an analytic manifold with entire Grauert tube in section 2, we found a polynomial $P(T)$ in \cref{eqn:Chens-poly-bd} that bounds the integral of the geodesic counting function $n_{T}(x,y)$ from above.

Denote by ${\rm Id}$ the identity function and define the constant $V_{M} := \int\limits_{M} {\rm Id_{M}}dx$.
Substituting all this information into Ma\~n\'e's formula in \cref{eqn:mane} yields:
\begin{eqnarray*}
\h &=& \lim\limits_{T \rightarrow \infty} \frac{1}{T} \log \int\limits_{M\times M}n_{T}(x,y)
\, dx\, dy  \\
  &\leq&   \lim\limits_{T \rightarrow \infty} \frac{1}{T}\log \int\limits_{M} P(T)
\, dx \\
  &=&  \lim\limits_{T \rightarrow \infty} \frac{1}{T} \log\, \left[ P(T)\cdot V_{M} \right] = 0 
\end{eqnarray*}

Observe that we could have also used the formula in \cref{eqn:htop-Jacobi-formula} above, to argue in a similar way using the polynomial bound given by \cref{eqn:Chens-poly-bd} on the growth of the Jacobi fields.
Therefore $\h =0$, as claimed. \end{proof}

\medskip

\begin{proof}[Proof of \cref{thm:3D-EntireGrauert}.]
 First observe that Anderson and Paternain \cite[Theorem B]{AndersonPaternain03}, under the assumption that the Geometrization Conjucture holds for $3$-manifolds, showed that such a $Y$ with $\h =0$ is modelled on ${\mathbb S}^3,\, {\mathbb S}^2\times {\mathbb E},\, {\mathbb E}^3$ or on the nilpotent Lie group $Nil^3$.
 As we now know the Geometrization Conjecture is true \cite{Hamilton82, Perelman02, Perelman03}, we may use this result directly.
 Milnor \cite[Theorem 2.4]{Milnor76} showed that a manifold modeled on $Nil^3$ has both positive and negative sectional curvatures.
This rules out $Nil^3$, because manifolds with entire Grauert tubes have nonnegative sectional curvature \cite{LempertSzoke91}.

 Now we shall see that every geometric manifold modelled on ${\mathbb S}^3, {\mathbb S}^2\times {\mathbb E}$, or ${\mathbb E}^3$ admits an entire Grauert tube.
 Sz\H{o}ke showed that compact rank one symmetric spaces admit entire Grauert tubes \cite{Szoke91}. This covers the case of geometric manifolds modelled on ${\mathbb S}^3$.
 Aguilar showed that taking quotients under isometries preserves the property of having entire Grauert tube \cite{Aguilar01}.
This implies the result for the ${\mathbb S}^2\times {\mathbb E}$, and ${\mathbb E}^3$ geometric manifolds.
 For the case of ${\mathbb S}^2\times {\mathbb E}$, every such geometric manifold is finitely covered by $S^2\times S^1$, where the result already holds for the product metric.
 Likewise, by Bieberbach's theorem every flat manifold modelled on ${\mathbb E}^3$ is finitely covered by a flat $3$-torus $T^3$. 
 Therefore, in both cases we may appeal to Aguilar's results to conclude that all geometric manifolds modelled on ${\mathbb S}^3, {\mathbb S}^2\times {\mathbb E}$, or ${\mathbb E}^3$ admit an entire Grauert tube.

 For the second assertion, we recall that Biswas and Mj \cite[Theorem 1.3]{BiswasMj15} proved that a closed $3$-manifold $Y$ admits a good complexification if and only if one of the following holds; $Y$ admits a flat metric, or $Y$ admits a metric of constant positive curvature, or $Y$ is covered by the (metric) product of a round $S^2$ and the real line. These are precisely the equivalences proved above for having an entire Grauert tube, corresponding in the same order to the geometric manifolds modelled on ${\mathbb E}^3$, ${\mathbb S}^3$, or ${\mathbb S}^2\times {\mathbb E}$.
\end{proof}

\medskip

\begin{proof} [Proof of \cref{thm:4D-simply-connected-EntireGrauert}.]
In all cases, by \cref{thm:entire-Grauert-null-htop} we know there exists a metric $g$ such that $\h=0$.
 \begin{enumerate}[(i)]
  \item Paternain showed that a smooth $4$-manifold $X$ with a metric such that $\h=0$ is homeomorphic one of the manifolds of this list \cite{Paternain00}.
    \item  Recall that Paternain and Petean produced a list of such manifolds with null entropy \cite[Theorem E]{PaternainPetean04}. 
  It remains to rule out $Nil^3\times {\mathbb E}$ from their list.
  Milnor showed that nilpotent manifolds have mixed curvatures \cite{Milnor76}.
  Therefore, by Lempert--Sz\H{o}ke's result that manifolds with entire Grauert tube have nonnegative sectional curvature \cite{LempertSzoke91} we may rule out geometric $Nil^3\times {\mathbb E}$ manifolds.
  \item This is precisely the conclusion of another result of Paternain and Petean \cite[Theorem F]{PaternainPetean04}, under the hypothesis that $X$ admits a metric $g$ with $\h =0$.
 \end{enumerate}
This concludes the proof. \end{proof}
\medskip

\begin{proof} [Proof of \cref{thm:4D-elliptic-surfaces-EntireGrauert}.]
    It is a theorem of Wall that such an $S$ has a geometric structure that is compatible with the complex structure \cite[Theorem 4.7]{Wall86}.
    Now we rely on the Enriques--Kodaira classification to separate the rest of the proof into several cases.
    When the Kodaira dimension $\kappa$ of $S$ equals $1$, the fundamental group of $S$ grows exponentially. 
    This is incompatible with \cref{cor:entire-Grauert-subexp-pi1}.
    The case of $\kappa =0$ and odd first Betti number $b_1$ corresponds to the geometry $Nil^3\times E$, which has both positive and negative sectional curvatures \cite{Milnor76}, and hence cannot have entire Grauert tube by Lempert--Sz\H{o}ke's well known result \cite{LempertSzoke91}.
    The remaining cases are,
    \begin{enumerate}[(i)]
        \item $\kappa =0$, $b_1$ even, corresponding to ${\mathbb E}^4$;
        \item $\kappa = -\infty$, $b_1$ even, corresponding to ${\mathbb S}^2\times {\mathbb E}^2$, and;
        \item $\kappa =-\infty$, $b_1$ odd, corresponding to ${\mathbb S}^3\times {\mathbb E}$.
    \end{enumerate}
    By Aguilar's constructions \cite{Aguilar01}, the metrics defined by each of these geometric structures, ${\mathbb S}^3\times {\mathbb E}$, ${\mathbb S}^2\times {\mathbb E}^2$, and ${\mathbb E}^4$, have entire Grauert tubes.
\end{proof}

\medskip

\begin{proof} [Proof of \cref{thm:infinite-pi1-EntireGrauert}.]
    By \cref{cor:entire-Grauert-subexp-pi1} the fundamental group of $X$ has subexponential growth type.
    Given the potential existence of finitely presented groups of intermediate growth, we assume that $\pi_1(X)$ has polynomial growth.
    Then, by a result of Paternain and Petean \cite[Theorem 4.5]{PaternainPetean06}, and \cref{thm:entire-Grauert-null-htop}, we conclude that $X$ is finitely covered by either a manifold homeomorphic to $S^3\times S^1$, or by a manifold homeomorphic to an infranilmanifold, or by a manifold $s$-cobordant to  $S^2\times T^2$.
    Hillman \cite{Hillman07} further showed that the case where the covering is $s$-cobordant to $S^2\times T^2$ can be improved to conclude that there exists a finite covering of $X$ homeomorphic to $S^2\times T^2$.
    The case of having a covering homeomorphic to an infranilmanifold implies that the manifold is a $K(\pi_1(X),1)$-space. 
    Recall that $X$ has nonnegative sectional curvature \cite{LempertSzoke91}.
    Therefore, by a result of Cheeger and Gromoll \cite[Theorem 9.5]{CheegerGromoll72}, $X$ must be flat.
\end{proof}

\medskip

\begin{proof} [Proof of \cref{thm:simply-conn-5mfds-EGT}.]
If $(M,g)$ has entire Grauert tube, then $\h =0$ by \cref{thm:entire-Grauert-null-htop}.
By a result of Paternain and Petean \cite[Theorem E]{PaternainPetean04} a closed simply connected $5$-manifold admits a smooth metric $g$ with $\h = 0$ if and only if it is diffeomorphic to $S^5$, $S^3\times S^2$, the nontrivial $S^3$-bundle over $S^2$ or the Wu manifold ${\rm SU}(3)/{\rm SO}(3)$.

The manifolds $S^5$, $S^3\times S^2$,  and ${\rm SU}(3)/{\rm SO}(3)$ all have both entire Grauert tubes and good complexifications \cite{PatrizioWong91, Aguilar01, Totaro03}.
Following Barden, let $X_{\infty}$ denote the nontrivial $S^3$-bundle over $S^2$ \cite{Barden65}.  
It was shown by Pavlov that $X_{\infty}$ is a biquotient \cite{Pavlov04} (see also the construction provided by DeVito \cite{DeVito11}). 
We thank Totaro for explaining the following description of $X_{\infty}$ as a biquotient, with $S^1$ acting on $S^3 \times S^3$. 
Let $u$ be in ${\rm SO}(2)\simeq S^1$, and
consider the free action of $S^1$ on the first $S^3$ factor defined by $\left( \begin{smallmatrix}
u & 0 \\
0 & u
\end{smallmatrix}\right)\in {\rm SO}(4)$.
Similarly, $S^1$ acts on the second $S^3$ factor by $u\mapsto \left( \begin{smallmatrix}
{\rm Id} & 0 \\
0 & u
\end{smallmatrix}\right)$.
As the first action is free, and its quotient is $S^2$, the projection onto the first coordinate gives the fibration $S^3 \hookrightarrow X_{\infty}\to S^2$. 
Therefore $X_{\infty} = (S^3 \times S^3)/S^1$ is a biquotient, and it has a good complexification by Totaro's result \cite{Totaro03}.

Let us briefly verify that each of the manifolds $S^5$, $S^3\times S^2$, the nontrivial $S^3$-bundle over $S^2$ and ${\rm SU}(3)/{\rm SO}(3)$ has an entire Grauert tube. 
Recall that all round spheres $S^{n}$ have entire Grauert tube, so this covers $S^{5}$. 
Aguilar observed that a product of two manifolds with entire Grauert tube also has entire Grauert tube \cite[\S 5.1\, (2)]{Aguilar01}. 
Hence $S^3\times S^2$ has entire Grauert tube. 
Let $G$ be a compact Lie group and $H$ a closed subgroup of $G$.
A bi-invariant metric $g$ on $G$ induces a metric $g_{n}$ on the homogeneous space $G/H$, called the normal metric.
A Riemannian symmetric space of nonnegative sectional curvature is a normal homogeneous Riemannian manifold.
It is well known that the Wu manifold ${\rm SU}(3)/{\rm SO}(3)$ is a symmetric space of non-negative curvature, so it is a normal homogeneous Reimannian manifold. 
Now we recall a result of Sz\H{o}ke, showing that  normal homogeneous Reimannian manifolds have entire Grauert tube \cite[Theorem 2.2]{Szoke98}.
Therefore ${\rm SU}(3)/{\rm SO}(3)$ has entire Grauert tube.
Observe that biquotients of normal homogeneous Reimannian  manifolds are normal homogeneous themselves.
Hence, the biquotient presentation of $X_{\infty}$ above implies it has entire Grauert tube, again by Sz\H{o}ke's  result \cite[Theorem 2.2]{Szoke98} (see also \cite[Corollary 4.2]{Aguilar01}).
\end{proof}

\providecommand{\bysame}{\leavevmode\hbox to3em{\hrulefill}\thinspace}

\end{document}